\documentclass[12pt]{amsart}
\usepackage[utf8]{inputenc}
\usepackage[T2A]{fontenc}
\usepackage{fullpage}
\usepackage[english]{babel}

\usepackage{amsmath}
\usepackage{amssymb}
\usepackage{amsthm}
\usepackage{amscd}
\usepackage{multirow}
\usepackage{bigdelim}

\theoremstyle{plain}
\newtheorem{Theorem}{Theorem}
\newtheorem{Proposition}{Proposition}
\newtheorem{lemma}{Lemma}

\newtheorem{Conjecture}{Conjecture}

\theoremstyle{definition}

\newtheorem{remark}{Remark}

\newtheorem{nsfactor}{}

\def \appendixref #1 {Appendix,~no.~\ref{#1}}

\def \prooflike #1 {{\sc {#1.}} }
\def \endproof {$\Box$}

\DeclareMathOperator{\rk}{rk}
\DeclareMathOperator{\diag}{diag}
\DeclareMathOperator{\Lie}{Lie}
\DeclareMathOperator{\GL}{GL}
\DeclareMathOperator{\SL}{SL}
\DeclareMathOperator{\SO}{SO}
\DeclareMathOperator{\OrthogonalGroup}{O}
\DeclareMathOperator{\Spin}{Spin}

\DeclareMathOperator{\Sp}{Sp}

\DeclareMathOperator{\GroupA}{A}
\DeclareMathOperator{\GroupB}{B}
\DeclareMathOperator{\GroupG}{G}
\DeclareMathOperator{\GroupF}{F}
\DeclareMathOperator{\GroupH}{H}
\DeclareMathOperator{\GroupT}{T}
\DeclareMathOperator{\GroupR}{R}

\DeclareMathOperator{\GroupP}{P}

\DeclareMathOperator{\Tangent}{T}

\DeclareMathOperator{\Mat}{Mat}

\def \isom {\cong}
\def \AA {\mathbb{A}}

\def \ZZ {\mathbb{Z}}

\def \kk {\Bbbk}

\def \catquot {/\hspace{-3pt}/}
\def \dcosets #1#2#3 {#1 \hspace{-1pt} \backslash\hspace{-3pt}\backslash\hspace{-0.8pt}{#2}\hspace{-1pt}\slash\hspace{-3pt}\slash #3 \hspace{1pt}}
\def \predcosets #1#2#3 {#1 \backslash{#2}\slash #3 \hspace{1pt}}
\def \oneptspace { \left\{ \rm pt \right\} }

\begin{document}

\subjclass[2000]{14L30,14M17}
\title{Spherical subgroups and double coset varieties}
\address{Department of Higher Algebra, Faculty of Mechanics and Mathematics, Lomonosov Moscow State University, Leninskie Gory 1, GSP-1, Moscow 119991, Russia}
\email{aanisimov@inbox.ru}
\author{Artem~B.~Anisimov}

\maketitle

\begin{abstract}
Let~$\GroupG$ be a connected reductive algebraic group,~$\GroupH \subsetneq \GroupG$ a reductive subgroup and~$\GroupT \subset \GroupG$ a maximal torus. It is well known that if charactersitic of the ground field is zero, then the homogeneous space~$\GroupG/\GroupH$ is a smooth affine variety, but never an affine space. The situation changes when one passes to double coset varieties~$\dcosets{\GroupF}{\GroupG}{\GroupH}$. In this paper we consider the case of~$\GroupG$ classical and~$\GroupH$ connected spherical and prove that either the double coset variety~$\dcosets{\GroupT}{\GroupG}{\GroupH}$ is singular, or it is an affine space. We also list all pairs~$\GroupH \subset \GroupG$ such that~$\dcosets{\GroupT}{\GroupG}{\GroupH}$ is an affine space.
\end{abstract}

\section{Introduction}

The construction of homogeneous spaces~$\GroupG / \GroupH$ has a natural generalisation: one can take another subgroup~$\GroupF \subset \GroupG$ and, instead of~$\GroupH$-cosets, consider~$(\GroupF,\GroupH)$-cosets, namely, sets~$\GroupF g \GroupH$. Double cosets in~$\GroupG$ play an important role in a wide variety of problems concerning actions of~$\GroupG$. For instance, in~\cite{RichardsonRohrleSteinberg} it was shown that if~$\GroupG$ is simple and $\GroupP = LQ \subset \GroupG$ is a parabolic subgroup with abelian unipotent radical~$Q$, $L$ being a Levi subgroup of~$\GroupP$, then the number of~$L$-orbits in~$Q$ is the same as cardinality of~$\predcosets{\GroupP}{\GroupG}{\GroupP}$, which is finite. Another example is enumeration of simple modules over simple groups~$\GroupG$ with finite number of~$\GroupG$-orbits on subspaces of dimension~$k$: if~$V$ is a simple~$\GroupG$-module satisfying this property and~$\GroupP_k$ is the stabiliser in~$\SL(V)$ of a~$k$-subspace of~$V$ then, as shown in~\cite{GLMS}, the set~$\predcosets{\GroupG}{\SL(V)}{\GroupP_k}$ is finite. Many other examples can be found in an expository paper~\cite{Seitz}.

In this paper we consider \textit{double coset varieties}. The variety~$\dcosets{\GroupF}{\GroupG}{\GroupH}$ is defined to be the categorical quotient of~$\GroupG$ with respect to the action of~$\GroupF \times \GroupH$ given by the formula~$(f,h) \circ g = fgh^{-1}$. We get interesting problems by posing the simplest questions: when~$\dcosets{\GroupF}{\GroupG}{\GroupH}$ exists and when the action~$\GroupF \times \GroupH : \GroupG$ is locally transitive?

If no restrictions are imposed upon~$\GroupF$ and~$\GroupH$ then the very existence of~$\dcosets{\GroupF}{\GroupG}{\GroupH}$ is not guaranteed; one can apply results from~\cite{ArzhCelik} to find out when this variety exists. We limit ourselves to the case of reductive subgroups~$\GroupF$ and~$\GroupH$. In this case the double coset variety exists and coincides with the spectrum of the algebra~${}^{\GroupF} \kk\left[ \GroupG \right] {}^{\GroupH}$ of functions on~$\GroupG$ invariant with respect to the described action of~$\GroupF \times \GroupH$. Remark that~$\dcosets{\GroupF}{\GroupG}{\GroupH}$ parametrises closed double $(\GroupF,\GroupH)$-cosets in~$\GroupG$. A result of Luna~\cite{LunaClosedOrbits} asserts that the action~$\GroupF \times \GroupH : \GroupG$ is stable, hence~$\dcosets{\GroupF}{\GroupG}{\GroupH}$ parametrises generic~$(\GroupF,\GroupH)$-cosets.

The second question, namely, existence of a dense coset~$\GroupF g \GroupH$, has also been extensively studied. In our setting, when~$\GroupF$ and~$\GroupH$ are reductive and when the ground field has zero characterstic, the question on density of~$\textrm{FH}$ reduces to the question on existence of the decomposition~$\GroupG = \textrm{FH}$. Indeed, since the action~$\GroupF \times \GroupH : \GroupG$ is stable, the dense coset~$\textrm{FH}$ is closed and coincides with~$\GroupG$. Compact simple Lie groups admitting a decomposition~$\GroupG = \GroupF\GroupH$ are classified in~\cite{OnishchikInclRel}. Subsequent paper~\cite{OnishchikDecomp} provides this classification in context of reductive Lie groups.

It is interesting to see what happens when one replaces actions of~$\{e\} \times \GroupH : \GroupG$ with actions of~$\GroupF \times \GroupH : \GroupG$~--- what are the properties of actions~$\GroupH : \GroupG$ that change and what are those ones that remain the same? For instance, all orbits of actions~$\GroupH : \GroupG$ are closed; according to~\cite{LunaClosedOrbits}, the actions of~$\GroupF \times \GroupH$ retain this property for generic orbits. But some of properties change radically. For example, as shown in~\cite{KraftPopov}, if groups~$\GroupH \subsetneq \GroupG$ are reductive, then the homogeneous space~$\GroupG / \GroupH$ is never an affine space\footnote{It is important that characteristic of the ground field be zero. For~${\rm char}\ \kk = 2$ there is an example of a transitive action~$\SL_2 : \mathbb{A}^2$, see~\cite{Knop}. This remark has been communicated to us by W. van der Kallen.}. Meanwhile, in~\cite{ArzhGayf} it was observed that~$\dcosets{\GroupT}{\SL_4}{\Sp_4}$, with~$\GroupT$ being a maximal torus of~$\SL_4$, is the affine plane. If one does not require the subgroups~$\GroupF$ and ~$\GroupH$ to be reductive then numerous examples of this kind can be constructed. Indeed, if~$\GroupF$ and~$\GroupH$ are \textit{excellent}, then, as proved in~\cite{VinbergGindikin}, the variety~$\dcosets{\GroupF}{\GroupG}{\GroupH}$ is an affine space. Recall that a spherical subgroup~$\GroupH \subseteq \GroupG$ is said to be excellent if the weight semigroup~$\Lambda_+\left( \GroupG / \GroupH \right)$ is generated by disjoint linear combinations of fundamental weights, that is, no fundamental weight appears with non-zero coefficient in two or more generators; the weight semigroup~$\Lambda_+\left( \GroupG / \GroupH \right)$ consists of highest weights of simple~$\GroupG$-modules having non-trivial~$\GroupH$-invariant vectors.

It is natural to pose a question to describe double coset varieties that are ``the simplest ones'', that is, those varieties~$\dcosets{\GroupF}{\GroupG}{\GroupH}$ that are affine spaces. The similar problem has already been resolved for many classes of linear representations of reductive groups, see~\cite[Section 8]{VinbergPopov}. For double coset varieties no such classification exists at the time.

In this paper we prove an easily verifiable necessary condition for~$\dcosets{\GroupF}{\GroupG}{\GroupH}$ to be an affine space. In one special case, namely, when~$\GroupG$ is a classical group, $\GroupH \subset \GroupG$ is a connected spherical reductive subgroup and~$\GroupF \subset \GroupG$ is a maximal torus, we enumerate all pairs~$\GroupH \subseteq \GroupG$ such that the algebra~${}^{\GroupF} \kk\left[ \GroupG \right] {}^{\GroupH}$ is free. By ``classical'' we mean the groups~$\SL_n$, $\SO_n$ and~$\Sp_{2n}$. The described class of subgroups~$\GroupF$ and~$\GroupH$ resembles the class considered in~\cite{PanyushevActionOfMaxTorus} which deals with the case where~$\GroupF$ is a maximal torus and~$\GroupH$ is the stabiliser of highest weight vector of a simple~$\GroupG$-module.

Now let us proceed to the main results of the paper.
\begin{Theorem}\label{criterion}
Let~$\GroupG$ be a classical algebraic group, $\GroupT \subset \GroupG$ be a maximal torus,~$\GroupH \subset \GroupG$ be a connected spherical reductive subgroup and let~$\pi : \GroupG \rightarrow \dcosets{\GroupT}{\GroupG}{\GroupH}$ be the quotient morphism. Then the double coset variety~$\dcosets{\GroupT}{\GroupG}{\GroupH}$ is an affine space if and only if the image~$\pi(e)$ of the identity element is a regular point.
\end{Theorem}

Theorem~\ref{criterion} resembles the criterion for linear representations of reductive groups to have free algebra of invariants~\cite[Section 4.4]{VinbergPopov}.

\begin{Theorem}\label{classification}
Let~$\GroupG$ be a classical algebraic group, $\GroupT \subset \GroupG$ be a maximal torus,~$\GroupH \subset \GroupG$ be a connected spherical reductive subgroup. Then the double coset variety~$\dcosets{\GroupT}{\GroupG}{\GroupH}$ is an affine space if and only if the groups~$\GroupG$ and~$\GroupH$ are listed in the following table:

\begin{center}\begin{tabular}{lll}
$\GroupG$ \phantom{MMMMMM}		&		$\GroupH$	\phantom{MMMMMMMM}	&		$\dim \dcosets{\GroupT}{\GroupG}{\GroupH}$	\\
\hline
$\SL_{n+1}$					&		${\rm S}\left( \GL_n \times \GL_1 \right)$		&		$n$																								\\
$\SL_4$							&		$\Sp_4$															&		$2$																								\\
$\SO_{2n+1}$				&		$\SO_{2n}$														&		$n$																								\\
$\SO_{2n}$						&		$\SO_{2n-1}$													&		$n-1$																							\\
$\SO_4$							&		$\GL_2$															&		$1$																								\\
\hline
$\SO_8$							&		$\Spin_7$														&		$3$																								\\
$\SO_6$							&		$\GL_3$															&		$3$																								\\
$\SO_4$							&		$\SO_2 \times \SO_2$										&		$2$																								\\
$\SO_3$							&		$\GL_1$															&		$1$																								\\
$\Sp_4$							&		$\Sp_2 \times \Sp_2$										&		$2$																								\\
\end{tabular}\end{center}
\end{Theorem}

\begin{remark}
In certain sense the bottom part of the table duplicates the top one. Indeed, the pair $\GL_3 \hookrightarrow \SO_6$ is an image of~${\rm S}\left( \GL_3 \times \GL_1 \right) \subset \SL_4$ under the two-fold covering $\SL_4 \rightarrow \SO_6$. The pair~$\Sp_2 \times \Sp_2 \subset \Sp_4$, on the contrary, is a two-fold covering of $\SO_4 \subset \SO_5$. Similarly, the cases~$\SO_2 \times \SO_2 \subset \SO_4$ and~$\GL_1 \subset \SO_3$ both reduce to~${\rm S}(\GL_1 \times \GL_1) \subset \SL_2$ and the case~$\Spin_7 \subset \SO_8$ reduces to~$\SO_7 \subset \SO_8$.
\end{remark}

Theorems~\ref{criterion} and~\ref{classification} are proved by walking through Krämer's list~\cite{Kramer} of connected spherical reductive subgroups in simple groups which consists of all symmetric pairs and~$12$ other pairs. We filter out those double coset varieties~$\dcosets{\GroupT}{\GroupG}{\GroupH}$ that are not affine spaces by pointing out their singular points. To this end we use the following proposition.
\begin{Proposition}\label{necessary-condition-of-smoothness}
Let~$\GroupF, \GroupH \subseteq \GroupG$ be reductive subgroups and~$\pi : \GroupG \rightarrow \dcosets{\GroupF}{\GroupG}{\GroupH}$ be the quotient morphism. Suppose that the double coset~$\GroupF\GroupH$ is closed in~$\GroupG$. Let~$Z$ be the categorical quotient for the action~$\GroupF \cap \GroupH : \Lie\GroupG / \left( \Lie\GroupF + \Lie\GroupH \right)$ induced by the adjoint action of~$\GroupF \cap \GroupH$ on~$\Lie\GroupG$. Then the point~$\pi(e) \in \dcosets{\GroupF}{\GroupG}{\GroupH}$ is regular if and only if~$Z$ is an affine space.
\end{Proposition}

Applying Proposition~\ref{necessary-condition-of-smoothness} we determine whether~$\pi(e)$ is a regular point or not. In cases when~$\pi(e)$ is regular we check that~$\dcosets{\GroupT}{\GroupG}{\GroupH}$ is an affine space.

D.~I.~Panyushev suggested the following observation: for all pairs~$\GroupH \subset \GroupG$ listed in Theorem~\ref{classification}, except~$\SO_2 \times \SO_2 \subset \SO_4$, the weight semigroup~$\Lambda_+\left( \GroupG / \GroupH \right)$ is generated by one element.
\begin{Conjecture}\label{rank-conjecture}
Let~$\GroupG$ be a simple algebraic group and $\GroupH \subset \GroupG$ a connected spherical reductive subgroup. If~$\dcosets{\GroupT}{\GroupG}{\GroupH}$ is an affine space then~$\rk \Lambda_+\left( \GroupG / \GroupH \right) = 1$.
\end{Conjecture}

If Conjecture~\ref{rank-conjecture} is true then Theorem~\ref{classification} lists all Krämer pairs that have free algebra~${}^{\GroupT} \kk[ \GroupG ]^{\GroupH}$, not only ones with~$\GroupG$ classical. See Section~3.5 for details. 

This paper does not consider Krämer's pairs with exceptional groups~$\GroupG$ because in this case applying Proposition~\ref{necessary-condition-of-smoothness} is technically more difficult. It seems that rather than carrying out calculations for exceptional groups, it would be preferable to give an à priori proof of Conjecture~\ref{rank-conjecture} and Theorem~\ref{criterion}.

Throughout this text the ground field~$\kk$ is supposed to be algebraically closed and of characteristic zero. All topological terms refer to Zariski topology. Simple roots and fundamental weights are numbered as in~\cite{VinbergOnischik}.

The author would like to thank I.~V.~Arzhantsev for stating the problem and helpful discussions. Valuable suggestions of D.~I.~Panyushev have considerably simplified some of originally employed proofs and improved the overall structure of the paper. R.~S.~Avdeev has also made numerous useful suggestions.

\section{Necessary condition of smoothness of double coset varieties}

Let us begin by proving that double coset varieties depend only on conjugacy classes of subgroups~$\GroupF$ and~$\GroupH$.
\begin{Proposition}\label{conjugation-invariance}
Let~$\GroupF, \GroupH \subseteq \GroupG$ be arbitrary closed subgroups of an algebraic group~$\GroupG$ and~$x, y \in \GroupG$ be two elements. Let~$\GroupF^\prime = x \GroupF x^{-1}$, $\GroupH^\prime = y \GroupH y^{-1}$. Then the algebras~${}^{\GroupF} \kk\left[\GroupG\right]{}^{\GroupH}$ and~${}^{\GroupF^\prime} \kk\left[\GroupG\right]{}^{\GroupH^\prime}$ are isomorphic. In particular, if~$\GroupF$ and~$\GroupH$ reductive, the varieties~$\dcosets{\GroupF}{\GroupG}{\GroupH}$ and~$\dcosets{\GroupF^\prime}{\GroupG}{\GroupH^\prime}$ are isomorphic.
\end{Proposition}
\begin{proof}
The isomorphism takes a function~$\varphi \in {}^{\GroupF} \kk\left[\GroupG\right]{}^{\GroupH}$ to~$\varphi^\prime (g) = \varphi \left( x^{-1}gy \right)$.
\end{proof}

The above proposition shows that with no loss of generality one may assume the subgroups~$\GroupF$ and~$\GroupH$ to \textit{have maximal intersection}, that is, for every~$g \in \GroupG$ dimension of~$g \GroupF g^{-1} \cap \GroupH$ is not greater than dimension of~$\GroupF \cap \GroupH$.

\begin{Proposition}\label{type-of-max-intersection}
Let~$\GroupH \subseteq \GroupG$ be a connected reductive subgroup, $\GroupF \subseteq \GroupG$ be a maximal torus (resp. a maximal unipotent subgroup, or a Borel subgroup). If~$\GroupH$ and~$\GroupF$ have maximal intersection, then~$\GroupH \cap \GroupF$ is a maximal torus (resp. a maximal unipotent subgroup or a Borel subgroup) of~$\GroupH$.
\end{Proposition}
\begin{proof}
Let~$\GroupF_0 \subseteq \GroupH$ be a subgroup of the same type as~$\GroupF$ (maximal torus, maximal unipotent subgroup or a Borel subgroup). There exists~$g \in \GroupG$ such that~$g \GroupF_0 g^{-1} \subseteq \GroupF$. For this element we have~$\GroupF_0 \subset g^{-1} \GroupF g \cap \GroupH$. To complete the proof we need to show that~$g^{-1} \GroupF g \cap \GroupH$ is connected.

If~$\GroupF$ is a maximal torus in~$\GroupG$, then~$g^{-1} \GroupF g \cap \GroupH$ is commutative and centralises~$\GroupF_0$. By~\cite[IX.24.1]{Humphreys}, it coincides with~$\GroupF_0$. Unipotent subgroups are automatically connected. Finally, if~$\GroupF$ is a Borel subgroup of~$\GroupG$ then~$g^{-1} \GroupF g \cap \GroupH$ is connected since Borel subgroups are maximal not only as connected solvable subgroups, but also as closed solvable subgroups, see~\cite[VIII.23.1]{Humphreys}.
\end{proof}

Proposition~\ref{necessary-condition-of-smoothness} requires the double coset~$\GroupF\GroupH$ to be closed. Let us show that this restriction is superfluous.
\begin{Proposition}\label{maximal-intersection-is-reductive}
Suppose that reductive subgroups~$\GroupF, \GroupH \subset \GroupG$ have maximal intersection. Then~$\GroupF\GroupH$ is closed in~$\GroupG$ and~$\GroupF \cap \GroupH$ is reductive.
\end{Proposition}
\begin{proof}
It is clear that the stabiliser of~$g \in \GroupG$ in~$\GroupF \times \GroupH$ is isomorphic to~$\GroupF \cap g \GroupH g^{-1}$. It follows that~$\dim\left( \GroupF g \GroupH \right) = \dim\GroupG - \dim\left( \GroupF \cap g \GroupH g^{-1} \right)$. Since the subgroups~$\GroupF$ and~$\GroupH$ have maximal intersection, the double coset~$\GroupF e \GroupH$ is an orbit of minimal dimension, therefore it is closed. By Matsushima's criterion~\cite[Theorem 4.17]{VinbergPopov}, the intersection~$\GroupF \cap \GroupH$ is reductive.
\end{proof}

\prooflike{Proof of Proposition~\ref{necessary-condition-of-smoothness}}
Let us first find the tangent space of~$\GroupF\GroupH$ at~$e$. Obviously,~$\Lie\GroupF, \Lie\GroupH \subseteq \Tangent_e \left( \GroupF\GroupH \right)$, hence we have~$\Lie\GroupF + \Lie\GroupH \subseteq \Tangent_e \left( \GroupF\GroupH \right)$. The double coset~$\GroupF\GroupH$ is an orbit of~$e$, hence~$\dim \Tangent_e(\GroupF\GroupH) = \dim(\GroupF\GroupH) = \dim(\GroupF \times \GroupH) - \dim(\GroupF \cap \GroupH)$. On the other hand,
\begin{align*}
\dim(\Lie\GroupF + \Lie\GroupH) =
\dim\Lie\GroupF + \dim\Lie\GroupH - \dim(\Lie\GroupF \cap \Lie\GroupH) =\\
= \dim(\GroupF \times \GroupH) - \dim(\GroupF \cap \GroupH).
\end{align*}
Thus,~$\Tangent_e(\GroupF\GroupH) = \Lie\GroupF + \Lie\GroupH$.

The double coset~$\GroupF\GroupH$ is a closed orbit of~$\GroupF \times \GroupH$, therefore, by Luna's theorem~\cite{LunaSlicesEtales} on slice étalé, we have the following commutative diagram:
\begin{equation*}
\begin{CD}
\left( \GroupF \times \GroupH \right) \ast_{\GroupR} S				@>>>			\GroupG\\
@VVpV																@VVV\\
S \catquot \GroupR							@>f>>		\GroupG \catquot \left( \GroupF \times \GroupH \right)		&	=		\dcosets{\GroupF}{\GroupG}{\GroupH}.
\end{CD}
\end{equation*}
In this diagram~$\GroupR$ denotes the stabiliser in~$\GroupF \times \GroupH$ of point~$e$, the slice~$S$ is an open neighbourhood of the origin in the slice module~$N$,~$p : N \rightarrow N \catquot \GroupR$ is the quotient morphism and the morphism~$f$ is étale.

For étale morphisms~$\varphi : X \rightarrow Y$ we have~$x \in X^{reg} \Leftrightarrow \varphi(x) \in Y^{reg}$. Thus, regularity of~$\pi(e)$ is equivalent to regularity of~$p(0)$. By~\cite[Proposition~4.11]{VinbergPopov} the latter is equivalent to the fact that~$\kk\left[ N \right]^{\GroupR}$ is free.

It remains to prove that the representation~$\GroupR : N$ is isomorphic to the representation~$\GroupF \cap \GroupH : \Lie\GroupG / \left( \Lie\GroupF + \Lie\GroupH \right)$ induced by the adjoint action of~$\GroupF \cap \GroupH$ on $\Lie\GroupG$. The tangent algebra~$\Lie\GroupG$ has the following decomposition into~$\GroupR$-modules
\begin{equation*}
\Lie\GroupG = \Tangent_e \left( \GroupF\GroupH \right) \oplus N = \left( \Lie\GroupF + \Lie\GroupH \right) \oplus N.
\end{equation*}
It is clear that~$\GroupR = \left\{ \left( f,f^{-1} \right)\ |\ f \in \GroupF \cap \GroupH \right\}$. Therefore the action~$\GroupR : \Lie\GroupG$ is isomorphic to the adjoint action of~$\GroupF \cap \GroupH$ on $\Lie\GroupG$, hence its restriction to~$N$ is isomorphic to the action~$\GroupF \cap \GroupH : \Lie\GroupG / \left( \Lie\GroupF + \Lie\GroupH \right)$.
\endproof

Later in our reasoning we will have to check whether some specific linear representations of tori have free algebras of invariants. To this end, let us introduce necessary notation and prove a monotonicity result concerning linear  representations of tori.

Let~$\GroupT$ be an algebraic torus acting in a vector space~$V$. Consider the weight decomposition~$V = \bigoplus_{i=1}^{n} V_{\chi_i}$. It will be convenient to suppose that all weight spaces~$V_{\chi_i}$ are one-dimensional, but weights~$\chi_i$ can have multiplicities. Denote~$A(\GroupT, V)$ \textit{a semigroup of linear relations between weights~$
\left\{ \chi_i \right\}$}:
\begin{equation*}
A\left( \GroupT, V \right) = \left\{ \left( a_1, \dots, a_n \right) \in \ZZ_{\geq 0}^n\ \big|\ \sum\limits_{i=1}^{N} a_i \chi_i
 = 0 \right\}.
\end{equation*}

Denote~$x_1, \dots, x_n$ coordinates in a basis of~$V$ that consists of~$\GroupT$-weight vectors. It is clear that elements of~$A$ are~$n$-tuples~$(a_1, \dots, a_n)$ such that monomials~$x_1^{a_1} \cdots x_n^{a_n}$ belong to~$\kk\left[ V \right]^{\GroupT}$. This observation renders the following proposition obvious.
\begin{Proposition}\label{abhaendigkeit}
The algebra~$\kk\left[ V \right]^{\GroupT}$ is free if and only if the semigroup~$A\left( \GroupT, V \right)$ is free.
\end{Proposition}

The following lemma will often be used to prove that certain representations~$\GroupT : V$ have singular quotients~$V \catquot \GroupT$.

\begin{lemma}\label{monotonicity-lemma}
Let~$\GroupT = \GroupT_0 \times \GroupT_1$ be an algebraic torus acting in a vector space~$V$ and let $\left\{ \chi_1, \dots, \chi_s, \chi_{s+1}, \dots, \chi_{s+r} \right\}$ be the weights of~$V$ with respect to this action. Suppose that the characters~$\chi_i$ with~$i \leq s$ have trivial restrictions to~$\GroupT_1$. Denote~$U = \bigoplus\nolimits_{i=1}^{s} V_{\chi_i}$. Under these assumptions, if the categorical quotient~$U \catquot \GroupT_0$ is singular, then so is~$V \catquot \GroupT$.
\end{lemma}
\begin{proof}
Since~$\chi_i$ with~$i \leq s$ have trivial restrictions to~$\GroupT_1$, we have~$U \catquot \GroupT = U \catquot \GroupT_0$, hence~$U \catquot \GroupT$ is singular; thus, every collection of generators of~$\kk[ U ]^{\GroupT}$ has algebraic relations. According to~\cite[Proposition~1.3]{PopovSyzigies}, every minimal set of generators of the algebra~$\kk[U]^{\GroupT}$ can be turned into a minimal set of generators of~$\kk[V]^{\GroupT}$ by adding several invariant functions; such minimal set of generators of~$\kk[V]^{\GroupT}$ also has algebraic relations. This implies that the quotient~$V \catquot \GroupT$ is not an affine space, hence it is singular.
\end{proof}

As we will see later on, numerous families of triples~$\GroupT, \GroupH \subseteq \GroupG$ have slice modules that satisfy conditions of Lemma~\ref{monotonicity-lemma}. For such series of triples Lemma~\ref{monotonicity-lemma} will be used to show that if a double coset variety~$\dcosets{\GroupT_r}{\GroupG_r}{\GroupH_r}$ is not an affine space when~$r = r_0$ then it is not an affine space for all~$r > r_0$.

\section{Proof of Theorems~\ref{criterion} and~\ref{classification}}

\

\subsection{Spherical subgroups of special linear groups}
We denote~$\GroupT_n \subset \SL_n$ the subgroup of the diagonal matrices. When there is no possibility for confusion we write~$\GroupT$ instead of~$\GroupT_n$. By maximal torus we always mean the maximal torus consisting of diagonal matrices. Propositions~\ref{conjugation-invariance} and~\ref{type-of-max-intersection} show that we can assume~$\GroupT_n \cap \GroupH_n$ to be maximal tori in groups~$\GroupH_n$. By singular algebra we mean an algebra of regular functions on a singular affine variety.

\begin{Proposition}\label{case-s-gl-gl-in-sl}
Let~$n$ and~$m$ be two positive integers,~$n \geq m$,~$\GroupG = \SL_{n+m}$ and $\GroupH_{n,m} = {\rm S}\left( \GL_n \times \GL_m \right) \subset \SL_{n+m}$. The algebra~${}^{\GroupT} \kk\left[ \GroupG \right]{}^{\GroupH_{n,m}}$ is free if~$m=1$ and singular otherwise.
\end{Proposition}
\begin{proof}
Take~$m \geq 2$. The slice module~$N_{n,m} = \Lie\SL_{n+m} / \left( \Lie\GroupT + \Lie\GroupH_{n,m} \right)$ has weights~$\pm( e_i - e_j )$ with~$i \leq n$, $j \leq m$. As shown in~\appendixref{appendix-s-gl-gl-in-sl} , the categorical quotient~$N_{n,m} \catquot T$ with~$m \geq 2$ is singular. By Proposition~\ref{necessary-condition-of-smoothness}, the algebra~${}^{\GroupT} \kk\left[ \GroupG \right] {}^{\GroupH}$ is singular when~$m \geq 2$.

Now let us show that if~$m=1$ then the algebra~${}^{\GroupT} \kk[\GroupG]^{\GroupH}$is free. Consider the group~$\GroupR \subset \GroupH$, which is~$\SL_n$ acting on first~$n$ basis vectors and consider its action on space~$\Mat_{n+1}$ of square matrices of size~$n+1$. An element~$r \in \GroupR$ transforms a matrix~$X \in \Mat_{n+1}$ according to the rule~$r \circ X = Xr^{-1}$; one easily checks that the algebra of invariants~$\kk\left[ \Mat_{n+1} \right]^{\SL_n}$ is generated by following functions:
\begin{itemize}
\item $M_i$~---~minor of~$X$ obtained by removing the~$i$-th row and the last column,
\item $y_i$~---~element of the last column~$X$ located in the~$i$-th row.
\end{itemize}

The group~$\SL_{n+1} \subset \Mat_{n+1}$ is an~$\GroupR$-invariant closed subset in~$\Mat_{n+1}$, therefore all functions in~$\kk\left[ \SL_{n+1} \right]^{\GroupR}$ are restrictions of~$\GroupR$-invariant functions on~$\Mat_{n+1}$. Thus,~$\kk\left[ \SL_{n+1} \right]^{\SL_n} = \kk\left[ M_i,\ y_j\ |\ i,j = 1, \dots, n+1 \right]$. The group~$\GroupH$ contains also all diagonal matrices in~$\SL_{n+1}$; by considering action of these elements we find that~$\kk\left[ \SL_{n+1} \right]^{\GroupH} = \kk\left[ M_iy_j\ |\ i,j = 1, \dots, n+1 \right]$.

The functions~$M_i$ and~$y_j$ are semiinvariant with respect to action of the diagonal torus~$\GroupT \subset \SL_{n+1}$ by left multiplications:
\begin{equation*}
\left\{\begin{array}{rcl}
M_i	&	\mapsto	&	t_i^{-1} M_i, \\
y_i	&	\mapsto	&	t_i y_i.
\end{array}\right.
\end{equation*}
Therefore~${}^{\GroupT} \kk\left[ \SL_{n+1} \right] {}^{\SL_n} = \kk\left[ M_iy_i\ |\ i = 1, \dots, n+1 \right]$. The listed generators satisfy one linear relation.
\begin{equation*}
\sum\limits_{k=1}^{n} \left( -1 \right)^{n+1+k} M_k y_k = \det X = 1.
\end{equation*}
This relation shows that~${}^{\GroupT} \kk\left[ \SL_{n+1} \right] {}^{\SL_n} = \kk\left[ M_iy_i\ |\ i = 1, \dots, n \right]$. It is clear that these~$n$ generators are algebraically independent, hence the algebra~${}^{\GroupT} \kk\left[ \SL_{n+1} \right] {}^{\SL_n}$ is free.
\end{proof}

\begin{Proposition}\label{case-sl-sl-in-sl}
Let~$n$ and~$m$ be two positive integers,~$n > m \geq 1$,~$\GroupG = \SL_{n+m}$ and $\GroupH = \SL_n \times \SL_m \subset \SL_{n+m}$. The algebra~${}^{\GroupT} \kk\left[\GroupG\right]{}^{\GroupH}$ is singular for all~$n$ and~$m$.
\end{Proposition}
\begin{proof}
A simple modification of the previous reasoning shows that for all~$m$ the minimal number of generators of the algebra~$\kk\left[ N \right]^{\GroupT}$ is greater than~$\dim N \catquot \GroupT$. Thus, all varieties~$\dcosets{\SL_n}{\SL_{n+m}}{\SL_m}$ are singular.
\end{proof}

\begin{Proposition}\label{case-sp-in-sl}
Let~$\GroupG = \SL_{2m}$ and $\GroupH = \Sp_{2m} \subset \SL_{2m}$. The algebra~${}^{\GroupT} \kk\left[\GroupG\right]{}^{\GroupH}$ is free if~$m=1$ or~$m=2$ and singular otherwise.
\end{Proposition}
\begin{proof}
The slice module~$\Lie\SL_{2m} / \left( \Lie\GroupT_{2m} + \Lie\Sp_{2m} \right)$ is
\begin{equation*}
N_{2m} = \left\{
	\left.\left(\begin{array}{c|c}
		\huge 0		&		B	\\
		\hline
		C				&		A
	\end{array}\right)
	\right|
	B = -B^s, C = -C^s, \mbox{ the diagonal of } A \mbox{ is zero}
\right\}.
\end{equation*}
In the above equality~$A^s$ denotes transposition of matrix~$A$ with respect to its secondary diagonal.

With respect to~$\GroupT_{2m} \cap \Sp_{2m}$ the block~$A$ has weights~$\varepsilon_i - \varepsilon_j$, where~$i \neq j$, the block~$B$ has weights~$\varepsilon_i + \varepsilon_j$ with~$j > i$ and the block~$C$ has weights~$-\varepsilon_i - \varepsilon_j$ with~$j > i$. Thus, the slice representations~$\GroupT_{2m} \cap \Sp_{2m} : N_{2m}$ satisfy the conditions of Lemma~\ref{monotonicity-lemma}.

According to~\appendixref{appendix-sp6-in-sl6} , the slice module for~$\Sp_6 \subset \SL_6$ has a singular categorical quotient~$N_6 \catquot (\GroupT_6 \cap \Sp_6)$, hence the algebra~${}^{\GroupT} \kk\left[ \SL_6 \right]^{\Sp_6}$ is singular. By Lemma~\ref{monotonicity-lemma}, the algebra~${}^{\GroupT_{2m}} \kk\left[ \SL_{2m} \right] {}^{\Sp_{2m}}$ is singular if~$m \geq 3$.

When~$m=1$ this assertion is trivial because~$\Sp_2 = \SL_2$ and~$\dcosets{\GroupT_2}{\SL_2}{\Sp_2} = \oneptspace$.

When~$m=2$ we begin by considering the action of the group~$\Sp_4$ on space of square matrices of order~$4$ given by the formula~$g \circ X = Xg^{-1}$. One easily checks that~$\kk\left[ \Mat_4 \right]^{\Sp_4} = \kk\left[ \left( u_i, u_j \right)\ |\ i,j = 1, \dots, 4 \right]$ where~$\left( u_i, u_j \right)$ denotes pairing of rows~$i$ and~$j$ of matrix~$X$ with respect to the bilinear form preserved by~$\Sp_4$. The listed generators are semiinvariant with respect to the action of the diagonal torus of~$\SL_4$ by left multiplications. We have
\begin{equation*}
{}^{\GroupT_4} \kk\left[ \SL_4 \right] {}^{\Sp_4} = \kk\left[ \left( u_i, u_j \right) \left( u_k, u_l \right)\ |\ \left\{i,j,k,l\right\} = \left\{1,2,3,4\right\} \right].
\end{equation*}
As we can see, the algebra~${}^{\GroupT_4} \kk\left[ \SL_4 \right] {}^{\Sp_4}$ is generated by three elements. These elements satisfy one linear relation.
\begin{equation*}
\left( u_1, u_2 \right)\left( u_3, u_4 \right) -
\left( u_1, u_3 \right)\left( u_2, u_4 \right) +
\left( u_1, u_4 \right)\left( u_2, u_3 \right) = \det X = 1.
\end{equation*}
Indeed, the expression on the left-hand side is a skew-symmetric bilinear function of rows and it equals~$1$ when~$u_i$ are rows of the identity matrix; hence it coincides with~$\det X$. This linear relation shows that one of the generators can be omitted. The remaining two generators are algebraically independent, hence the algebra~${}^{\GroupT_4} \kk\left[ \SL_4 \right] {}^{\Sp_4}$ is free.
\end{proof}

\begin{Proposition}\label{case-k-sp-in-sl}
Let~$\GroupG = \SL_{2m+1}$ and~$\GroupH_{2m} = \kk^\times \cdot \Sp_{2m} \subset \SL_{2m+1}$. Then the algebra ${}^{\GroupT_{2m+1}} \kk\left[\GroupG\right]{}^{\GroupH_{2m}}$ is free if~$m=1$ and singular otherwise.
\end{Proposition}
\begin{proof}
Maximal torus of group~$\GroupH_{2m}$ consists of the diagonal matrices
\begin{equation*}
\diag\left( tt_1, tt_2, \dots, tt_m, tt_m^{-1}, \dots, tt_1^{-1}, t^{-2m} \right).
\end{equation*}
It is clear that this torus is isomorphic to~$\left( \GroupT_{2m} \cap \Sp_{2m} \right) \times \kk^\times / \left\{ \pm 1 \right\}$. Let us regard the slice module~$N_{2m}$ as a module not over~$\GroupT_{2m+1} \cap \GroupH_{2m}$, but over its two-fold covering torus~$\widetilde{\GroupT}_{2m} = \left( \GroupT_{2m+1} \cap \Sp_{2m} \right) \times \kk^\times$. We regard characters~$\varepsilon_i$ of the torus~$\GroupT_{2m+1} \cap \Sp_{2m}$ as characters of~$\widetilde{\GroupT}_{2m}$ and denote~$\varepsilon$ the basis character which corresponds to the factor~$\kk^\times$. It is clear that the torus~$\widetilde{\GroupT}_{2m}$ acts in the slice module~$N_{2m}$ with the following weights:
\begin{equation*}
\begin{array}{ll}
\pm\left( \varepsilon_i \pm \varepsilon_j \right),		&		i < j, \\
\pm\left( (2m+1)\varepsilon \pm \varepsilon_i \right).
\end{array}
\end{equation*}
As we can see, the slice module~$\widehat{N}_{2m}$ corresponding to~$\Sp_{2m} \subset \SL_{2m}$ and the slice module~$N_{2m}$ satisfy the conditions of Lemma~\ref{monotonicity-lemma}. It follows from proof of the previous proposition that the categorical quotient~$N_{2m} \catquot \widetilde{\GroupT}_{2m}$ is singular if~$m \geq 3$. Two categorical quotients~$N_{2m} \catquot \widetilde{\GroupT}_{2m}$ and~$N_{2m} \catquot \left( \GroupT_{2m+1} \cap \GroupH_{2m} \right)$ coincide, hence the algebra~${}^{\GroupT_{2m+1}} \kk\left[ \SL_{2m+1} \right] {}^{\GroupH_{2m}}$ is singular if~$m \geq 3$.

When~$m=2$ we get a module considered in~\appendixref{appendix-k-sp-in-sl} , hence the algebra~${}^{\GroupT_{5}} \kk\left[ \SL_{5} \right] {}^{\GroupH_{4}}$ is singular.

The remaining case~$m=1$ has already been considered. Indeed, the group $\kk^{\times} \cdot \Sp_2 \subset \SL_3$ coincides with~${\rm S}\left( \GL_2 \times \GL_1 \right) \subset \SL_3$. By Proposition~\ref{case-s-gl-gl-in-sl}, the algebra~${}^{\GroupT_3} \kk\left[ \SL_3 \right] {}^{\GroupH_2}$ is free.
\end{proof}

\begin{Proposition}\label{case-so2m-in-sl2m}
Let~$\GroupG = \SL_{2m}$ and $\GroupH = \SO_{2m} \subset \SL_{2m}$. The algebra~${}^{\GroupT} \kk\left[\GroupG\right]{}^{\GroupH}$ is free if~$m=1$ and singular otherwise.
\end{Proposition}
\begin{proof}
The slice module~$\Lie\SL_{2m} / \left( \Lie\GroupT_{2m} + \Lie\SO_{2m} \right)$ is
\begin{equation*}
N_{2m} = \left\{
	\left.\left(\begin{array}{c|c}
		\huge 0		&		B	\\
		\hline
		C				&		A
	\end{array}\right)
	\right|
	B = B^s, C = C^s, \mbox{ the diagonal of } A \mbox{ is zero}
\right\}.
\end{equation*}

Thus, the slice module~$\widehat{N}_{2m}$ corresponding to~$\Sp_{2m} \subset \SL_{2m}$ and the slice module~$N_{2m}$ are isomorphic. From the proof of Proposition~\ref{case-sp-in-sl} and from Lemma~\ref{monotonicity-lemma} it follows that~${}^{\GroupT_{2m}} \kk\left[ \SL_{2m} \right] {}^{\SO_{2m}}$ is singular if~$m \geq 3$.

If~$m=2$ then the slice module~$N_4$ has weights~$\pm \varepsilon_1 \pm \varepsilon_2$, $\pm 2\varepsilon_1$ and~$\pm 2\varepsilon_2$. According to~\appendixref{appendix-so4-in-sl4} , the categorical quotient~$N_4 \catquot \GroupT$ is singular, hence the algebra~${}^{\GroupT_4} \kk\left[ \SL_4 \right] {}^{\SO_4}$ is singular.

Finally, let us consider the case~$\SO_2 \subset \SL_2$. The group~$\SO_2$ coincides with the diagonal torus of~$\SL_2$ and hence it coincides with~${\rm S}\left( \GL_1 \times \GL_1 \right) \subset \SL_2$. By Proposition~\ref{case-s-gl-gl-in-sl}, the algebra~${}^{\GroupT_2} \kk\left[ \SL_2 \right] {}^{\SO_2}$ is free. In this case the double coset variety is~$\AA^1$.
\end{proof}

\begin{Proposition}
Let~$\GroupG = \SL_{2m+1}$ and $\GroupH = \SO_{2m+1} \subset \SL_{2m+1}$. The algebra~${}^{\GroupT} \kk\left[\GroupG\right]{}^{\GroupH}$ is singular for all~$m$.
\end{Proposition}
\begin{proof}
The slice module~$\widehat{N}_{2m}$ corresponding to~$\SO_{2m} \subset \SL_{2m}$ and the slice module~$N_{2m+1} = \Lie\SO_{2m+1} / \left( \Lie\GroupT_{2m+1} + \Lie\SO_{2m+1} \right)$ satisfy the conditions of Lemma~\ref{monotonicity-lemma}. Hence, from the proof of Proposition~\ref{case-so2m-in-sl2m} it follows that the algebra~${}^{\GroupT_{2m+1}} \kk\left[ \SL_{2m+1} \right] {}^{\SO_{2m+1}}$ is singular if~$m \geq 2$.

If~$m=1$ then the slice module~$N_3$ has weights~$\pm \varepsilon_1$ and~$\pm 2\varepsilon_1$. According to~\appendixref{appendix-so3-in-sl3} , the quotient~$N_3 \catquot (\GroupT_3 \cap \SO_3)$ is singular, hence the algebra~${}^{\GroupT_3} \kk\left[ \SL_3 \right]{}^{\SO_3}$ is singular.
\end{proof}

\begin{Proposition}
Let~$\GroupG = \SL_{2m+1}$ and $\GroupH = \Sp_{2m} \subset \SL_{2m+1}$. The algebra~${}^{\GroupT} \kk\left[\GroupG\right]{}^{\GroupH}$ is singular for all~$m$.
\end{Proposition}
\begin{proof}
The slice module~$N = \Lie\SL_{2m+1} / \left( \GroupT + \Lie\Sp_{2m} \right)$ is
\begin{equation*}
\left(
\begin{array}{ccc|cc}
	\multicolumn{1}{c}{\multirow{2}*{\Large 0}}	& \multicolumn{1}{c|}{}		&
	\multicolumn{1}{c|}{\multirow{2}*{\Large B}}													&		x_1\\
	&	\multicolumn{1}{c|}{}		&																				&		\vdots \\
\cline{1-3}
	\multicolumn{1}{c}{\multirow{2}*{\Large C}}	& \multicolumn{1}{c|}{}		&
	\multicolumn{1}{c|}{\multirow{2}*{\Large A}}													&		\vdots\\
	&	\multicolumn{1}{c|}{}		&																				&		x_{2m} \\
\hline
y_1 \dots		&&		\dots		y_{2m}	&		0
\end{array}
\right),
\end{equation*}
where~$A$, $B$ and $C$ are square matrices of order~$m$; $B = B^s$, $C = C^s$ and~$A$ has only zero elements in its primary diagonal. Therefore, the maximal torus of~$\Sp_{2m}$ acts in the slice module with weights~$\pm \varepsilon_i \pm \varepsilon_j$ and $\pm \varepsilon_k$ with~$i,j,k \leq m$. As we can see, the slice modules corresponding to~$\Sp_{2m} \subset \SL_{2m+1}$ satisfy the conditions of Lemma~\ref{monotonicity-lemma}, thus it suffices to prove that~$\dcosets{\GroupT}{\SL_3}{\Sp_2}$ is singular.

The case~$m=1$ has already been considered: it is~$\SL_2 \subset \SL_3$. As we know from proof of Proposition~\ref{case-sl-sl-in-sl}, the categorical quotient~$N \catquot \GroupT$ is singular. Thus, for every~$m \geq 1$ the variety~$\dcosets{\GroupT}{\SL_{2m+1}}{\Sp_{2m}}$ is singular.
\end{proof}

\

\subsection{Spherical subgroups of symplectic groups} Let us settle notations and definitions concerning symplectic groups. In vector space~$\kk^{2n}$ we fix a basis denoted $e_1, \dots, e_n, e^\prime_n, \dots, e^\prime_1$. We define the symplectic group~$\Sp_{2n}$ as the group of isometries of a skew-symmetric bilinear form~$\omega$ such that~$\omega\left( e_i, e^\prime_i \right) = -\omega\left( e^\prime_i, e_i \right) = 1$ and all other pairings of basis vectors are zero. The described choice of symplectic groups simplifies calculations because when defined in this way, the groups~$\Sp_{2m}$ have maximal tori that consist of the diagonal matrices~$\diag\left( t_1, \dots, t_n, t_n^{-1}, \dots, t_1^{-1} \right)$.

When considering groups~$\Sp_{2n} \times \Sp_{2m} \subset \Sp_{2n+2m}$, we use another convention and define~$\Sp_{2n+2m}$ as isometry groups of orthogonal direct sum~$\kk^{2n} \oplus \kk^{2m}$ where both summands are equipped with aforementioned skew-symmetric bilinear forms.

\begin{Proposition}\label{case-sp-sp-in-sp}
Let~$n$ and~$m$ be two positive integers,~$\GroupG = \Sp_{2n+2m}$ and $\GroupH = \Sp_{2n} \times \Sp_{2m} \subset \Sp_{2n+2m}$. The algebra~${}^{\GroupT} \kk\left[\GroupG\right]{}^{\GroupH}$ is free if~$n=m=1$ and singular otherwise.
\end{Proposition}
\begin{proof}
The maximal torus~$\GroupT$ acts in~$\Lie\GroupG / \Lie\GroupH$ with weights~$\pm \varepsilon_i \pm \varepsilon_j$ with~$i \leq n$ and~$j > n$. So, the slice module corresponding to~$\Sp_{2n} \times \Sp_{2m} \subset \Sp_{2n+2m}$ and the one corresponding to~$\Sp_{2r} \times \Sp_{2s} \subset \Sp_{2r+2s}$ satisfy the conditions of Lemma~\ref{monotonicity-lemma} if~$r \geq n$ and~$s \geq m$.

In the case~$\Sp_2 \times \Sp_4 \subset \Sp_6$  the slice module~$N$ has weights~$\pm \varepsilon_1 \pm \varepsilon_2$ and~$\pm \varepsilon_1 \pm \varepsilon_3$. According to~\appendixref{appendix-sp2-sp4-in-sp6} , the categorical quotient~$N \catquot \GroupT$ is singular. Thus, the algebra~${}^{\GroupT} \kk\left[ \Sp_{2n+2m} \right] {}^{\Sp_{2n} \times \Sp_{2m}}$ is singular for all~$m$ and~$n$ except~$m=n=1$.

Let us now consider the remaining case~$\Sp_2 \times \Sp_2 \subset \Sp_4$. There is an isomorphism $\Sp_4 \isom \Spin_5$, hence we a have a two-fold covering~$\Sp_4 \rightarrow \SO_5$; this covering takes the subgroup~$\Sp_2 \times \Sp_2 \subset \Sp_4$ to~$(\Sp_2 \times \Sp_2) / \{\pm E\} \isom \SO_4$. Thus,~${}^{\GroupT} \kk[\Sp_4]^{\Sp_2 \times \Sp_2} = {}^{\GroupT} \kk[\SO_5]^{\SO_4}$. It will be shown in Proposition~\ref{case-so-2n-in-so-2n+1} that the latter algebra is free.
\end{proof}

\begin{Proposition}
Let~$\GroupG = \Sp_{2n+2}$ and $\GroupH = \Sp_{2n} \times \kk^\times \subset \Sp_{2n+2}$. The algebra~${}^{\GroupT} \kk\left[\GroupG\right]{}^{\GroupH}$ is singular for all~$n$.
\end{Proposition}
\begin{proof}
In this case it is convenient to define the group~$\Sp_{2n+2}$ as the group of isometries of the orthogonal direct sum~$\kk^{2n} \oplus \kk^2$ with summands equipped with standard skew-symmetric forms. When~$\Sp_{2n+2}$ is defined this way, the central torus~$\kk^\times$ of the group~$\GroupH$ acts as the diagonal torus~$\left\{ \diag\left( t, t^{-1} \right) \right\}$ on the summand~$\kk^2$.

The slice module~$N=\Lie\Sp_{2n+2} / \Lie\GroupH$ has weights~$\pm \varepsilon_i \pm \varepsilon_{n+1}$, $i \leq n$ and~$\pm 2\varepsilon_{n+1}$. According to~\appendixref{appendix-k-sp2n-in-sp2n2} , the categorical quotient~$N \catquot \GroupT$ is singular, hence the algebra~${}^{\GroupT} \kk\left[ \Sp_{2n+2} \right] {}^{\Sp_{2n} \times \kk^\times}$ is singular.
\end{proof}

\begin{Proposition}\label{case-gl-in-sp}
Let~$\GroupG = \Sp_{2n}$ and $\GroupH = \GL_n \hookrightarrow \Sp_{2n}$. The algebra~${}^{\GroupT} \kk\left[\GroupG\right]{}^{\GroupH}$ is free if~$n=1$ and singular otherwise.
\end{Proposition}
\begin{proof}
The group~$\GroupH=\GL_n$ imbeds into~$\Sp_{2n}$ as
\begin{equation*}
\GroupH = \left\{\left.
	\left(\begin{array}{c|c}
		\phantom{-}A\phantom{-}		&		0\\
		\hline
		0		&		\left( A^s \right)^{-1}
	\end{array}\right)
\ \right|\ A \in \GL_n\right\}.
\end{equation*}
Therefore the slice module~$N$ is
\begin{equation*}
\Lie\Sp_{2n} / \Lie\GroupH  = \left\{\left.
	\left(\begin{array}{c|c}
		0		&		B\\
		\hline
		C		&		0
	\end{array}\right)
\ \right|\ B^s = B,\ C^s = C \right\}.
\end{equation*}
Thus, the slice modules corresponding to~$\GL_n \hookrightarrow \Sp_{2n}$ satisfy the conditions of Lemma~\ref{monotonicity-lemma}.

In the case~$\GL_2 \hookrightarrow \Sp_4$ the slice module has weights~$\varepsilon_1 + \varepsilon_2$, $-\varepsilon_1 - \varepsilon_2$ and~$\pm 2\varepsilon_i$ where $i=1,2$. By~\appendixref{appendix-gl2-in-sp4} \ and Lemma~\ref{monotonicity-lemma}, the algebra~${}^{\GroupT} \kk\left[ \Sp_{2n} \right] {}^{\GL_n}$ is singular if~$n \geq 2$.

If~$n=1$ then the group~$\Sp_2$ coincides with~$\SL_2$ and the group~$\GroupH$ coincides with~${\rm S}\left( \GL_1 \times \GL_1 \right)$. By Proposition~\ref{case-s-gl-gl-in-sl}, the algebra~${}^{\GroupT} \kk\left[ \Sp_2 \right] {}^{\GL_1}$ is free.
\end{proof}

\

\subsection{Spherical subgroups of orthogonal groups} Let us settle notations and definitions concerning orthogonal groups. In even-dimensional vector space~$\kk^{2n}$ we fix a basis denoted $e_1, \dots, e_n, e^\prime_n, \dots, e^\prime_1$. We define the orthogonal group~$\OrthogonalGroup_{2n}$ to be the group of isometries of a symmetric bilinear form such that~$\left( e_i, e^\prime_i \right) = \left( e^\prime_i, e_i \right) = 1$ and all other pairings of basis vectors are zero. In odd-dimensional vector spaces~$\kk^{2n+1}$ we fix a basis~$e_1, \dots, e_n, e, e^\prime_n, \dots, e^\prime_1$ and consider the following symmetric bilinear form: pairings of~$e_i$ and~$e^\prime_i$ are the same as above, vector~$e$ is orthogonal to all other basis vectors and~$(e,e) = 1$. The described choice of orthogonal groups simplifies calculations because when defined in this way, the groups~$\SO_n$ have maximal tori that consist of the diagonal matrices~$\diag\left( t_1, \dots, t_n, 1, t_n^{-1}, \dots, t_1^{-1} \right)$. The unit component in the middle is present only in torus of~$\SO_{2n+1}$.

As with symplectic groups, when considering groups~$\SO_n \times \SO_m \subset \SO_{n+m}$, we define groups~$\OrthogonalGroup_{n+m}$ as isometry groups of orthogonal direct sum~$\kk^n \oplus \kk^m$ where both summands are equipped with aforementioned symmetric bilinear forms.

\begin{Proposition}\label{case-gl-in-so2n}
Let~$\GroupG = \SO_{2n}$ and $\GroupH = \GL_n \hookrightarrow \SO_{2n}$. The algebra~${}^{\GroupT} \kk\left[\GroupG\right]{}^{\GroupH}$ is free if~$n \leq 3$ and singular otherwise.
\end{Proposition}
\begin{proof}
The embedding~$\GL_n \hookrightarrow \SO_{2n}$ is the same as~$\GL_n \hookrightarrow \Sp_{2n}$, hence the slice module~$N$ is
\begin{equation*}
\Lie\SO_{2n} / \Lie\GroupH  = \left\{\left.
	\left(\begin{array}{c|c}
		0		&		B\\
		\hline
		C		&		0
	\end{array}\right)
\ \right|\ B^s = -B,\ C^s = -C \right\}.
\end{equation*}

The slice module~$N$ has weights~$\varepsilon_i + \varepsilon_j$ and $-\varepsilon_i - \varepsilon_j$, where~$i < j$, so the slice modules corresponding to~$\GL_n \hookrightarrow \SO_{2n}$ satisfy the conditions of Lemma~\ref{monotonicity-lemma}.

When~$n=4$ we get a module considered in~\appendixref{appendix-gl4-in-so8} . Thus the algebra~${}^{\GroupT_n} \kk\left[ \SO_{2n} \right] {}^{\GL_n}$ is singular if~$n \geq 4$.

The case~$n=3$ has already been considered. Indeed, one has the isomorphism~$\SO_6 \isom \SL_4 / \{ \pm E \}$ which takes~$\GroupH$ to~${\rm S}\left( \GL_3 \times \GL_1 \right) / \{ \pm E \}$. Remark that the generators of the algebra~${}^{\GroupT} \kk\left[ \SL_4 \right] {}^{{\rm S}\left(\GL_3 \times \GL_1\right)}$ that have been found in proof of Proposition~\ref{case-s-gl-gl-in-sl} are invariant with respect to multiplication by~$\pm E$, hence they they are, in fact, functions on~$\SO_6$. This shows that the algebra~${}^{\GroupT}\kk\left[ \SO_6 \right] {}^{\GL_3}$ is isomorphic to~${}^{\GroupT} \kk\left[ \SL_4 \right] {}^{{\rm S}\left(\GL_3 \times \GL_1\right)}$ and therefore free.

If~$n=2$ then we begin by applying a properly selected transformation from~$\GL_2$ that takes a generic matrix from~$\SO_4$ to a matrix
\begin{equation*}
\begin{pmatrix}
1		&		0		&		p				&		0		\\
0		&		1		&		0				&		-p		\\
u		&		0		&		1+up		&		0		\\
0		&		-u		&		0				&		1+up
\end{pmatrix}.
\end{equation*}
Further reducing an open subset of generic matrices we can assume that~$u \neq 0$. Applying a properly selected transformation from~$\GroupT \times \GL_2$, namely, the conjugation with~$\diag\left( u, 1, u^{-1}, 1 \right)$ we can replace a generic matrix with the following one:
\begin{equation*}
\begin{pmatrix}
1		&		0		&		up			&		0		\\
0		&		1		&		0				&		-up	\\
1		&		0		&		1+up		&		0		\\
0		&		-1		&		0				&		1+up
\end{pmatrix}.
\end{equation*}
It is clear that the function~$up$ extends to a~$\GroupT \times \GL_2$-invariant function~$F$ on the group~$\SO_4$, hence~${}^{\GroupT}\kk\left[ \SO_4 \right] {}^{\GL_2} = \kk\left[ F \right]$. This shows that this algebra is free.

The last remaining case~$n=1$ is trivial because~$\SO_2 = \left\{ \diag\left( t, t^{-1} \right)\ |\ t \in \kk^\times \right\}$ and we have~$\dcosets{\GroupT}{\SO_2}{\GL_1} = \oneptspace$.
\end{proof}

\begin{Proposition}
Let~$\GroupG = \SO_{2n}$ and $\GroupH = \SL_n \hookrightarrow \SO_{2n}$. The algebra~${}^{\GroupT} \kk\left[\GroupG\right]{}^{\GroupH}$ is singular for all~$n$.
\end{Proposition}
\begin{proof}
As in the case~$\GL_n \hookrightarrow \SO_{2n}$, the torus~$\GroupT_n \cap \SL_n$ acts in the slice module~$N_n$ with weights~$\varepsilon_i + \varepsilon_j$ and~$-\varepsilon_i - \varepsilon_j$ with~$i < j \leq n$, but this time the vectors~$\varepsilon_i$ are linearly dependent:~$\varepsilon_1 + \dots + \varepsilon_n = 0$. Because of this linear relation the slice modules corresponding to different pairs~$\SL_n \hookrightarrow \SO_{2n}$ do not satisfy the conditions of Lemma~\ref{monotonicity-lemma}. That is why we need to give a direct proof for every~$n$.

Let us point out some of linear relations between generators of semigroups~$A(N, \GroupT)$ that have to be included into minimal sets of their generators. Obviously, one has to include relations $\left( \varepsilon_i + \varepsilon_j \right) + \left( -\varepsilon_i - \varepsilon_j \right) = 0$ into minimal set of generators because these are the generators with minimal number of non-zero coefficients. Also there are two families of relations that have only coefficients~$0$ and~$1$ and that have to be included into minimal set of generators. These families are constructed in the following way: begin with the sum $\left( \varepsilon_1 + \varepsilon_{i_1} \right) + \left( \varepsilon_1 + \varepsilon_{j_1} \right)$ and add summands~$\left( \varepsilon_{i_1} + \varepsilon_{i_2} \right) + \left( \varepsilon_{j_1} + \varepsilon_{j_2} \right)$, $\left( \varepsilon_{i_2} + \varepsilon_{i_3} \right) + \left( \varepsilon_{j_2} + \varepsilon_{j_3} \right)$ and so on; if~$n$ is odd we add a final summand~$\left( \varepsilon_{i_s} + \varepsilon_{j_s} \right)$. The sum that is constructed in this way equals~$2\varepsilon_1 + \dots + 2\varepsilon_n = 0$. Similarly one constructs a family of relations~$\left( -\varepsilon_1 - \varepsilon_i \right) + \left( -\varepsilon_1 - \varepsilon_j \right) + \dots = 0$. The number of generators of~$A(N,\GroupT)$ that we have listed is greater than~$\dim N \catquot \GroupT$; thus, by Proposition~\ref{abhaendigkeit}, the algebra~${}^{\GroupT} \kk[\GroupG] {}^{\GroupH}$ is singular.
\end{proof}

\begin{Proposition}
Let~$\GroupG = \SO_{2n+1}$ and $\GroupH = \GL_n \hookrightarrow \SO_{2n+1}$. The algebra~${}^{\GroupT} \kk\left[\GroupG\right]{}^{\GroupH}$ is free if~$n=1$ and singular otherwise.
\end{Proposition}
\begin{proof}
Slice modules corresponding to~$\GL_n \hookrightarrow \SO_{2n+1}$ satisfy the conditions of Lemma~\ref{monotonicity-lemma}. When~$n=2$, a trivial modification of proof of Proposition~\ref{case-gl-in-sp} shows that the algebra~${}^{\GroupT} \kk\left[ \SO_{5} \right] {}^{\GL_2}$ is singular; thus, varieties~$\dcosets{\GroupT}{\SO_{2n+1}}{\GL_n}$ are singular if~$n \geq 2$.

If~$n=1$ we have~$\SO_3 \isom \SL_2 / \{ \pm E \}$ and this isomorphism takes~$\GroupH$ to~${\rm S}\left( \GL_1 \times \GL_1 \right) / \{ \pm E \}$. By Proposition~\ref{case-s-gl-gl-in-sl}, the algebra~${}^{\GroupT} \kk\left[ \SO_3 \right] {}^{\GL_1}$ is free.
\end{proof}

\begin{Proposition}\label{case-so-so-in-so}
Let~$n$ and~$m$ be two positive integers,~$n,m \geq 2$,~$\GroupG = \SO_{n+m}$ and $\GroupH = \SO_n \times \SO_m \subset \SO_{n+m}$. The algebra~${}^{\GroupT} \kk\left[\GroupG\right]{}^{\GroupH}$ is free if~$m=n=2$ and singular otherwise.
\end{Proposition}
\begin{proof}
The slice module~$\Lie\GroupG / (\Lie\GroupH + \Lie\GroupT)$ is
\begin{equation*}
N = \left\{ \left. \begin{pmatrix}0&B\\ C&0\end{pmatrix} \right|\ C = -B^s \right\},
\end{equation*}
where~$B$ is a~$n \times m$-matrix and~$C$ is a~$m \times n$-matrix.

Denote~$\varepsilon_i$ and~$\delta_j$ the standard basis characters of maximal tori of groups~$\SO_n$ and~$\SO_m$ respectively; we regard~$\varepsilon_i$ and~$\delta_j$ as characters of the maximal torus of~$\SO_{n+m}$.

Let both~$n$ and~$m$ be even and put~$n=2r$ and~$m=2s$. In this case the slice module has weights~$\pm \varepsilon_i \pm \delta_j$ where~$i=1, \dots, r$, $j=1, \dots, s$ and indices~$i$ and~$j$ can be equal. Note that this description of weights is also true when~$r=1$ or~$s=1$. As we can see, if~$r^\prime \geq r$ and~$s^\prime \geq s$ then the slice module corresponding to~$\SO_{2r} \times \SO_{2s} \subset \SO_{2r+2s}$ and the one corresponding to~$\SO_{2r^\prime} \times \SO_{2s^\prime} \subset \SO_{2r^\prime + 2s^\prime}$ satisfy the conditions of Lemma~\ref{monotonicity-lemma}. When we put~$n=2$ and~$m=4$ we obtain a module considered in~\appendixref{appendix-so2-so4-in-so6} . Therefore the algebra~${}^{\GroupT} \kk\left[ \SO_{2r + 2s} \right] {}^{\SO_{2r} \times \SO_{2s}}$ is singular if one of numbers~$r$ or~$s$ is greater than~$1$.

Let us consider the case~$\SO_2 \times \SO_2 \subset \SO_4$. The group~$\SO_2 \times \SO_2$ is the diagonal torus~$\GroupT$ of the group~$\SO_4$. Note that~$\SO_4 \isom \left( \SL_2 \times \SL_2 \right) / \{ \pm E \}$ and recall that the double coset variety~$\dcosets{\GroupT_2}{\SL_2}{\GroupT_2}$ is the affine line~$\AA^1$. Remark also that~$\GroupT_2 \times \GroupT_2$-invariant functions on~$\SL_2$ are invariant with respect to multiplication by~$\pm E$ and hence define functions on~$\SO_4$. Therefore we have
\begin{equation*}
\dcosets{\GroupT}{\SO_4}{\GroupT} = \dcosets{\left( \GroupT_2 \times \GroupT_2 \right)}{\left( \SL_2 \times \SL_2 \right)}{\left( \GroupT_2 \times \GroupT_2 \right)},
\end{equation*}
and the variety on the right-hand side is the affine plane~$\AA^2$.

Now consider the case where~$n = 2r$ is even, and~$m = 2s+1$ is odd. In this case the slice module has weights~$\pm \varepsilon_i \pm \delta_j$ and~$\pm \varepsilon_i$ where~$i=1, \dots, r$, $j=1, \dots, s$, hence the slice modules corresponding to different pairs~$\SO_{2r} \times \SO_{2s+1} \subset \SO_{2r+2s+1}$ satisfy the conditions of Lemma~\ref{monotonicity-lemma}. When~$r=s=1$ we get a module considered in~\appendixref{appendix-so2-so3-in-so5} . Therefore the varieties~$\dcosets{\GroupT}{\SO_{2r+2s+1}}{\SO_{2r} \times \SO_{2s+1}}$ are singular for all~$r,s$.

The last case when both~$n$ and~$m$ are odd is similar to the previous one. As in the previous case, all varieties~$\dcosets{\GroupT}{\SO_{2r+2s+2}}{\SO_{2r+1} \times \SO_{2s+1}}$ are singular.
\end{proof}

\begin{Proposition}\label{case-so-2n-1-in-so-2n}
Let~$\GroupG = \SO_{2n}$ and~$\GroupH = \SO_{2n-1}$. The algebra~${}^{\GroupT} \kk\left[ \SO_{2n} \right] {}^{\SO_{2n-1}}$ is free for all~$n$.
\end{Proposition}
\begin{proof}
Consider the standard basis~$e_1, \dots e_n, e^\prime_n, \dots, e^\prime_1$ in~$V = \kk^{2n}$ and let~$x_1, \dots, x_n$, $x^\prime_n, \dots, x^\prime_1$ be coordinates in this basis. The homogeneous space~$\SO_{2n} / \SO_{2n-1}$ is a quadric~$Z \subset \kk^{2n}$, $Z = \left\{ x_1x^\prime_1 + \dots + x_nx^\prime_n = 1 \right\}$. Obviously, the algebra~$\kk\left[ V \right]^{\GroupT}$ is freely generated by~$z_i = x_i x^\prime_i$. Thus,~$Z \catquot \GroupT \subset V \catquot \GroupT$ is a hyperplane $\{z_1 + \dots + z_n = 1\}$, hence~$\dcosets{\GroupT}{\SO_{2n}}{\SO_{2n-1}} = Z \catquot \GroupT$ is isomorphic to an affine space.
\end{proof}

\begin{Proposition}\label{case-so-2n-in-so-2n+1}
Let~$\GroupG = \SO_{2n+1}$ and~$\GroupH = \SO_{2n}$. The algebra~${}^{\GroupT} \kk\left[ \SO_{2n+1} \right] {}^{\SO_{2n}}$ is free for all~$n$.
\end{Proposition}
\begin{proof}
Take the standard basis $e_1, \dots e_n, e, e^\prime_n, \dots, e^\prime_1$ in~$V = \kk^{2n+1}$, let $x_1, \dots, x_n,  x, x^\prime_n, \dots, x^\prime_1$ be coordinates in this basis. As in the previous proof, $\SO_{2n+1} / \SO_{2n}$ is a quadric~$Z \subset \kk^{2n+1}$, this time~$Z = \left\{ x_1x^\prime_1 + \dots + x_nx^\prime_n + x^2 = 1 \right\}$. Since $\kk\left[ V \right]^{\GroupT} = \kk\left[ x_1x^\prime_1, \dots, x_nx^\prime_n, x \right]$, the subset~$Z \catquot \GroupT \subset V \catquot \GroupT$ is a cylinder over a parabola, hence it is isomorphic to an affine space.
\end{proof}

\begin{Proposition}\label{case-spin-in-so}
Let~$\GroupG = \SO_8$ and $\GroupH = \Spin_7$. The algebra~${}^{\GroupT} \kk\left[\GroupG\right]{}^{\GroupH}$ is free.
\end{Proposition}
\begin{proof}
All~$\Spin_8$-modules that are not~$\SO_8$-modules have no~$\GroupT$-invariant non-zero vectors, hence~${}^{\GroupT} \kk[\SO_8]^{\Spin_7} = {}^{\GroupT} \kk[\Spin_8]^{\Spin_7}$. Clearly, the latter algebra is isomorphic to~${}^{\GroupT} \kk[\SO_8]^{\SO_7}$, and, by Proposition~\ref{case-so-2n-1-in-so-2n}, the algebra~${}^{\GroupT} \kk[\SO_8]^{\SO_7}$ is free.
\end{proof}

\begin{Proposition}
Let~$\GroupG = \SO_9$ and $\GroupH = \Spin_7$. The algebra~${}^{\GroupT} \kk\left[\GroupG\right]{}^{\GroupH}$ is singular.
\end{Proposition}
\begin{proof}
The group~$\Spin_7$ embeds into the group~$\SO_9$ via $\Spin_7 \hookrightarrow \SO_8 \hookrightarrow \SO_9$, hence the slice module corresponding to~$\Spin_7 \hookrightarrow \SO_9$ is the one considered in~\appendixref{appendix-spin7-in-so9} . This shows that the algebra~${}^{\GroupT} \kk\left[ \SO_9 \right] {}^{\Spin_7}$ is singular.
\end{proof}

\begin{Proposition}
Let~$\GroupG = \SO_{10}$ and $\GroupH = \Spin_7 \times \SO_2$. The algebra~${}^{\GroupT} \kk\left[\GroupG\right]{}^{\GroupH}$ is singular.
\end{Proposition}
\begin{proof}
In this case the slice module coincides with a module considered in~\appendixref{appendix-spin7-so2-in-so10} . Thus the algebra~${}^{\GroupT} \kk\left[ \SO_{10} \right] {}^{\GroupH}$ is singular.
\end{proof}

\begin{Proposition}\label{case-g2-in-so7}
Let~$\GroupG = \SO_7$ or~$\GroupG = \SO_8$ and let $\GroupH = \GroupG_2$. The algebra~${}^{\GroupT} \kk\left[\GroupG\right]{}^{\GroupH}$ is singular.
\end{Proposition}
\begin{proof}
The group~$\GroupG_2$ embeds into~$\SO_8$ via $\GroupG_2 \hookrightarrow \SO_7 \hookrightarrow \SO_8$, hence the slice modules for~$\GroupG_2 \hookrightarrow \SO_7$ and~$\GroupG_2 \hookrightarrow \SO_8$ satisfy the conditions of Lemma~\ref{monotonicity-lemma}. The slice module for~$\GroupG_2 \hookrightarrow \SL_7$ is listed in~\appendixref{appendix-g2-in-so7} , hence ~$\dcosets{\GroupT}{\SO_7}{\GroupG_2}$ is singular. By Lemma~\ref{monotonicity-lemma}, $\dcosets{\GroupT}{\SO_8}{\GroupG_2}$ is also singular.
\end{proof}

\

\subsection{Proof of Theorems~\ref{criterion} and~\ref{classification}}\label{deducing-theorems}
We walk through Krämer's list of connected spherical reductive subgroups in simple groups. Every item of the list corresponds to one of Propositions~\ref{case-s-gl-gl-in-sl} through \ref{case-g2-in-so7}; each proposition first rejects pairs~$\GroupH \subset \GroupG$ that have singular algebras~${}^{\GroupT} \kk[ \GroupG ]^{\GroupH}$. To this end, Proposition~\ref{necessary-condition-of-smoothness} is applied; thus, singularity of algebra~${}^{\GroupT} \kk[ \GroupG ]^{\GroupH}$ follows from singularity of the point~$\pi(e) \in \dcosets{\GroupT}{\GroupG}{\GroupH}$. For the remaining pairs the point~$\pi(e)$ is regular and it turns out that in these cases~${}^{\GroupT} \kk[ \GroupG ]^{\GroupH}$ is free. In this way we obtain the list of Theorem~\ref{classification} and show that regularity of~$\pi(e)$ implies that~$\dcosets{\GroupT}{\GroupG}{\GroupH}$ is an affine space, thus proving Theorem~\ref{criterion}.

\

\subsection{A remark on exceptional groups~$\GroupG$}\label{section-exceptional-groups}
Assuming that Conjecture~\ref{rank-conjecture} is true we can assert that Theorem~\ref{classification} lists all Krämer pairs with simple~$\GroupG$ that have free algebra~${}^{\GroupT} \kk[ \GroupG ]^{\GroupH}$. Indeed, there are only two pairs with exceptional~$\GroupG$ and \hbox{$\rk \Lambda_+ (\GroupG / \GroupH) = 1$}, they are~$\GroupA_2 \subset \GroupG_2$ and~$\GroupB_4 \subset \GroupF_4$; let us check them.

\begin{Proposition}
Let~$\GroupG = \GroupF_4$ and~$\GroupH = \GroupB_4$. The algebra~${}^{\GroupT} \kk[\GroupG]^{\GroupH}$ is singular.
\end{Proposition}
\begin{proof}
The weights of the slice module~$N = \Lie\GroupG / \Lie\GroupH$ are~$(\pm\varepsilon_1 \pm\varepsilon_2 \pm\varepsilon_3 \pm\varepsilon_4)/2$. By~\appendixref{appendix-b4-in-f4} , the quotient~$N \catquot \GroupT$ is singular, hence~$\dcosets{\GroupT}{\GroupF_4}{\GroupB_4}$ is singular.
\end{proof}

\begin{Proposition}
Let~$\GroupG = \GroupG_2$ and~$\GroupH = \GroupA_2$. The algebra~${}^{\GroupT} \kk[\GroupG]^{\GroupH}$ is singular.
\end{Proposition}
\begin{proof}
The weights of the slice module~$N = \Lie\GroupG / \Lie\GroupH$ are~$\pm e_i$. According to \appendixref{appendix-a2-in-g2} , the quotient~$N \catquot \GroupT$ is singular, hence~$\dcosets{\GroupT}{\GroupG_2}{\GroupA_2}$ is singular.
\end{proof}

\section{Appendix. Some torus modules with singular quotients}
This appendix lists several linear representations of tori that arise as slice modules corresponding to pairs~$\GroupH \subset \GroupG$ considered in section 3. All of the listed representations have singular categorical quotients~$V \catquot \GroupT$. All cases are handled by a uniform reasoning, namely, we list~$\GroupT$-invariant monomials that have to be included into minimal sets of generators of algebras~$\kk[V]^{\GroupT}$ and in all cases it turns out that their number is greater than~$\dim V \catquot \GroupT$. We denote~$X_{\lambda}$ coordinates in~$V$ in a basis that consists of~$\GroupT$-weight vectors; so, for~$t \in \GroupT$ we have~$t \circ X_{\lambda} = \lambda(t)^{-1} X_{\lambda}$. In cases~\ref{appendix-s-gl-gl-in-sl},~\ref{appendix-g2-in-so7} and~\ref{appendix-a2-in-g2} the weights are linear combinations of vectors~$e_i$ which have one linear relation~$e_1 + \dots + e_r = 0$. If~$\GroupT$ is represented as a product of two subtori then we denote~$\varepsilon_i$ and~$\delta_j$ the basis characters of these subtori and regard~$\varepsilon_i$ and~$\delta_j$ as characters of~$\GroupT$.

Every item in the following list describes a linear representation~$\GroupT : V$ by enumerating the weights of~$V$ and then lists elements of the algebra~$\kk[V]^{\GroupT}$ that have to be included into its minimal set of generators.

\

\begin{nsfactor}\label{appendix-s-gl-gl-in-sl}
Torus of rank~$n+m-1$ ($n \geq m \geq 2$) acts with weights~$\pm( e_i - e_j )$ where $1 \leq i \leq n$ and $n+1 \leq j \leq n+m$:
\end{nsfactor}
\noindent $X_{e_i - e_j} X_{-e_i + e_j}$,\ 
$X_{e_i - e_j} X_{e_k - e_l} X_{-e_i + e_l} X_{-e_k + e_j}$,\ 
$X_{e_i - e_l} X_{e_k - e_j} X_{-e_i + e_j} X_{-e_k + e_l}$,\\
where~$1 \leq i,j \leq n$ and~$n+1 \leq k,l \leq n+m$.
\begin{nsfactor}\label{appendix-sp6-in-sl6}
Torus of rank~$3$ acts with weights~$\pm\varepsilon_i \pm\varepsilon_j$ where~$i < j$:
\end{nsfactor}
\noindent There are~$6$ invariant monomials of degree~$2$; they are $X_{\varepsilon_i + \varepsilon_j} X_{-\varepsilon_i - \varepsilon_j}$, $X_{\varepsilon_i - \varepsilon_j} X_{-\varepsilon_i + \varepsilon_j}$. Let us point out six monomials of degree~$3$: $X_{\varepsilon_1 + \varepsilon_2} X_{\varepsilon_3 - \varepsilon_1} X_{-\varepsilon_2 - \varepsilon_3}$, $X_{\varepsilon_1 + \varepsilon_2} X_{\varepsilon_3 - \varepsilon_2} X_{-\varepsilon_1 - \varepsilon_3}$ and four monomials obtained by cyclic permutations of indices~$1,2,3$.
\begin{nsfactor}\label{appendix-k-sp-in-sl}
Torus of rank~$3$ acts with weights $\pm\varepsilon_1 \pm\varepsilon_2$ and~$\pm\left( 5\varepsilon_3 \pm \varepsilon_i \right)$, where~$i=1,2$:
\end{nsfactor}
\noindent $X_{\varepsilon_1 - \varepsilon_2} X_{\varepsilon_2 - \varepsilon_1}$,\ 
$X_{\varepsilon_1 + \varepsilon_2} X_{-\varepsilon_1 - \varepsilon_2}$,\ 
$X_{5\varepsilon_3 \pm \varepsilon_i} X_{-5\varepsilon_3 \mp \varepsilon_i}$ and
$X_{5\varepsilon_3 \pm \varepsilon_i} X_{-5\varepsilon_3 \pm \varepsilon_j} X_{\mp \varepsilon_i \mp \varepsilon_j}$, where~$i \neq j$.
\begin{nsfactor}\label{appendix-so4-in-sl4}
Torus of rank~$2$ acts with weights $\pm\varepsilon_1 \pm\varepsilon_2$ and $\pm 2\varepsilon_i$:
\end{nsfactor}
\noindent $X_{\varepsilon_1 \pm \varepsilon_2} X_{-\varepsilon_1 \mp \varepsilon_2}$,\ 
$X_{2\varepsilon_i} X_{-2\varepsilon_i}$,\ 
$X_{\varepsilon_1 - \varepsilon_2} X_{-\varepsilon_1 - \varepsilon_2} X_{2\varepsilon_2}$,\ 
$X_{\varepsilon_2 - \varepsilon_1} X_{\varepsilon_1 + \varepsilon_2} X_{-2\varepsilon_2}$\ and
$X_{2\varepsilon_1} X_{2\varepsilon_2} X_{-\varepsilon_1 - \varepsilon_2}^2$.
\begin{nsfactor}\label{appendix-so3-in-sl3}
Torus of rank~$1$ acts with weights $\pm \varepsilon_1$ and~$\pm 2\varepsilon_1$:
\end{nsfactor}
\noindent $X_{\varepsilon_1} X_{-\varepsilon_1}$,\ 
$X_{2\varepsilon_1} X_{-2\varepsilon_1}$,\ 
$X_{\varepsilon_1}^2 X_{-2\varepsilon_1}$,\ 
$X_{-\varepsilon_1}^2 X_{2\varepsilon_1}$.
\begin{nsfactor}\label{appendix-sp2-sp4-in-sp6}
Torus of rank~$3$ acts with weights $\pm\varepsilon_1 \pm\varepsilon_2$ and~$\pm\varepsilon_1 \pm\varepsilon_3$:
\end{nsfactor}
\noindent $X_{\varepsilon_1 \pm \varepsilon_i} + X_{-\varepsilon_1 \mp \varepsilon_i}$,~$i=2,3$ and
$X_{\varepsilon_1 + \varepsilon_2} X_{\varepsilon_1 - \varepsilon_2} X_{-\varepsilon_1 + \varepsilon_3} X_{-\varepsilon_1 - \varepsilon_3}$,
$X_{\varepsilon_1 + \varepsilon_3} X_{\varepsilon_1 - \varepsilon_3} X_{-\varepsilon_1 + \varepsilon_2} X_{-\varepsilon - \varepsilon_2}$.
\begin{nsfactor}\label{appendix-k-sp2n-in-sp2n2}
Torus of rank~$n+1$ acts with weights $\pm\varepsilon_i \pm\varepsilon_{n+1}$ and~$\pm 2\varepsilon_{n+1}$, where~$i \leq n$:
\end{nsfactor}
\noindent $X_{\varepsilon_i \pm \varepsilon_{n+1}} X_{-\varepsilon_i \mp \varepsilon_{n+1}}$,\ 
$X_{\varepsilon_i + \varepsilon_{n+1}} X_{-\varepsilon_i + \varepsilon_{n+1}} X_{-2\varepsilon_{n+1}}$,\ 
$X_{\varepsilon_i - \varepsilon_{n+1}} X_{-\varepsilon_i - \varepsilon_{n+1}} X_{2\varepsilon_{n+1}}	$\ and 
$X_{2\varepsilon_{n+1}} X_{-2\varepsilon_{n+1}}$.
\begin{nsfactor}\label{appendix-gl4-in-so8}
Torus of rank~$4$ acts with weights $\pm\left( \varepsilon_i + \varepsilon_j \right)$, where~$i < j$:
\end{nsfactor}
\noindent $X_{\varepsilon_i + \varepsilon_j} X_{-\varepsilon_i - \varepsilon_j}$,\ 
$X_{\varepsilon_1 + \varepsilon_2} X_{\varepsilon_3 + \varepsilon_4} X_{-\varepsilon_1 - \varepsilon_3} X_{-\varepsilon_2 - \varepsilon_4}$,\ 
$X_{\varepsilon_1 + \varepsilon_2} X_{\varepsilon_3 + \varepsilon_4} X_{-\varepsilon_1 - \varepsilon_4} X_{-\varepsilon_2 - \varepsilon_3}$,\\
$X_{\varepsilon_1 + \varepsilon_3} X_{\varepsilon_2 + \varepsilon_4} X_{-\varepsilon_1 - \varepsilon_2} X_{-\varepsilon_3 - \varepsilon_4}$,\ 
$X_{\varepsilon_1 + \varepsilon_3} X_{\varepsilon_2 + \varepsilon_4} X_{-\varepsilon_1 - \varepsilon_4} X_{-\varepsilon_2 - \varepsilon_3}$.
\begin{nsfactor}\label{appendix-gl2-in-sp4}
Torus of rank~$2$ acts with weights $\pm\left( \varepsilon_1 + \varepsilon_2 \right)$ and~$\pm 2\varepsilon_i$:
\end{nsfactor}
\noindent $X_{\varepsilon_1 + \varepsilon_2} X_{-\varepsilon_1 - \varepsilon_2}$,\ 
$X_{2\varepsilon_i} X_{-2\varepsilon_i}$,\ 
$X_{\varepsilon_1 + \varepsilon_2}^2 X_{-2\varepsilon_1} X_{-2\varepsilon_2}$,\ 
$X_{-\varepsilon_1 - \varepsilon_2}^2 X_{2\varepsilon_1} X_{2\varepsilon_2}$.
\begin{nsfactor}\label{appendix-so2-so4-in-so6}
Product of tori of ranks~$1$ and~$2$ acts with weights $\pm \varepsilon_1 \pm \delta_j$:
\end{nsfactor}
\noindent $X_{\varepsilon_1 \pm \delta_i} X_{-\varepsilon_1 \mp \delta_i}$,\ 
$X_{\varepsilon_1 + \delta_1} X_{\varepsilon_1 - \delta_1} 	X_{-\varepsilon_1 + \delta_2} X_{-\varepsilon_1 - \delta_2}$,\ 
$X_{\varepsilon_1 + \delta_2} X_{\varepsilon_1 - \delta_2} 	X_{-\varepsilon_1 + \delta_1} X_{-\varepsilon_1 - \delta_1}$.
\begin{nsfactor}\label{appendix-so2-so3-in-so5}
Product of two tori of ranks~$1$ acts with weights $\pm \varepsilon_1 \pm \delta_1$ and $\pm \varepsilon_1$:
\end{nsfactor}
\noindent $X_{\varepsilon_1} X_{-\varepsilon_1}$,\ 
$X_{\varepsilon_1 \pm \delta_1} X_{-\varepsilon_1 \mp \delta_1}$,\ 
$X_{\varepsilon_1 + \delta_1} X_{\varepsilon_1 - \delta_1} X_{-\varepsilon_1}^2$,\ 
$X_{-\varepsilon_1 + \delta_1} X_{-\varepsilon_1 - \delta_1} X_{\varepsilon_1}^2$.
\begin{nsfactor}\label{appendix-spin7-in-so9}
Torus of rank~$3$ acts with weights $\pm 2\varepsilon_k$ and $\pm\varepsilon_1 \pm\varepsilon_2 \pm\varepsilon_3$:
\end{nsfactor}
\noindent $X_{2\varepsilon_k} X_{-2\varepsilon_k}$,\ 
$X_{\pm\varepsilon_1 \pm\varepsilon_2 \pm\varepsilon_3} X_{\mp\varepsilon_1 \mp\varepsilon_2 \mp\varepsilon_3}$,\ 
$X_{\varepsilon_1 \pm\varepsilon_2 \pm\varepsilon_3} X_{\varepsilon_1 \mp\varepsilon_2 \mp\varepsilon_3} X_{-2\varepsilon_1}	$.
\begin{nsfactor}\label{appendix-spin7-so2-in-so10}
Product of tori of ranks~$3$ and~$1$ acts with weights $\pm\varepsilon_1 \pm\varepsilon_2 \pm\varepsilon_3 \pm\delta_1$ and $\pm 2\varepsilon_k$:
\end{nsfactor}
\noindent $X_{2\varepsilon_k} X_{-2\varepsilon_k}$,\ 
$X_{\pm\varepsilon_1 \pm\varepsilon_2 \pm\varepsilon_3 + \delta} X_{\mp\varepsilon_1 \mp\varepsilon_2 \mp\varepsilon_3 - \delta}$,\ 
$X_{\varepsilon_1 \pm\varepsilon_2 \pm\varepsilon_3 \pm\delta} X_{\varepsilon_1 \mp\varepsilon_2 \mp\varepsilon_3 \mp\delta} X_{-2\varepsilon_1}	$.
\begin{nsfactor}\label{appendix-g2-in-so7}
Torus of rank~$2$ acts with weights $\pm\left( e_i + e_j \right)$, where~$i < j$:
\end{nsfactor}
\noindent $X_{e_i + e_j} X_{-e_i - e_j}$,\ 
$X_{e_1 + e_2} X_{e_1 + e_3} X_{e_2 + e_3}$,\ 
$X_{-e_1 - e_2} X_{-e_1 - e_3} X_{-e_2 - e_3}$.
\begin{nsfactor}\label{appendix-b4-in-f4}
Torus of rank~$4$ acts with weights $(\pm\varepsilon_1 \pm\varepsilon_2 \pm\varepsilon_3 \pm\varepsilon_4)/2$:
\end{nsfactor}
\noindent
$X_{\varepsilon_1 \pm\varepsilon_2 \pm\varepsilon_3 \pm\varepsilon_3} X_{-\varepsilon_1 \mp\varepsilon_2 \mp\varepsilon_3 \mp\varepsilon_3}$ and
$X_{\varepsilon_1 + \varepsilon_2 + \varepsilon_3 - \varepsilon_4} X_{\varepsilon_1 + \varepsilon_2 - \varepsilon_3 + \varepsilon_4} X_{-\varepsilon_1 - \varepsilon_2 + \varepsilon_3 + \varepsilon_4} X_{-\varepsilon_1 - \varepsilon_2 - \varepsilon_3 - \varepsilon_4}$.
\begin{nsfactor}\label{appendix-a2-in-g2}
Torus of rank~$2$ acts with weights~$\pm e_1, \pm e_2, \pm e_3$:
\end{nsfactor}
\noindent
$X_{e_i} X_{-e_i}$, $X_{e_1} X_{e_2} X_{e_3}$, $X_{-e_1} X_{-e_2} X_{-e_3}$.

\end{document}